\documentclass[11pt, a4paper, oneside]{smfart}

\usepackage[top=3cm, bottom=3cm, left=3cm, right=3cm]{geometry}
\usepackage{graphicx}
\usepackage[parfill]{parskip}
\usepackage{amsmath, amsthm, amssymb,amsfonts, mathrsfs}
\usepackage{wasysym}
\usepackage[utf8]{inputenc}
\usepackage[english, frenchb]{babel}
\usepackage[T1]{fontenc}
\usepackage{url}
\usepackage{verbatim}
\usepackage{indentfirst}
\usepackage{xypic}
\usepackage{braket}
\usepackage{fancyhdr}
\usepackage{enumerate}
\usepackage{textcomp}
\usepackage{stmaryrd}
\usepackage{xcolor}
\usepackage{colortbl}
\usepackage[colorlinks=true, linktoc=page, citecolor=blue, linkcolor=blue, urlcolor=blue]{hyperref}
\usepackage{fancyhdr}
\usepackage{smfthm}

\usepackage{garamondlibre}

\usepackage{setspace}
\NumberTheoremsAs{theos}

\usepackage[
backend=bibtex,
style=authoryear, 
citestyle=authoryear-comp,
maxnames=7,
sortcites=false %%%% to keep the order in \cite command
]{biblatex}
\usepackage{csquotes}

\DefineBibliographyExtras{french}{\restorecommand\mkbibnamefamily}

\DeclareDelimFormat{nameyeardelim}{\addcomma\space}

\DeclareNameAlias{sortname}{given-family}

\renewcommand{\bibnamedash}{\leavevmode\raise3pt\hbox to3em{\hrulefill}\space}

\AtEveryBibitem{%
  \clearfield{issn} % Remove issn
  \clearfield{isbn} % Remove isbn
  \clearfield{doi} % Remove doi
  \clearlist{language} %%% remove language
  \ifentrytype{online}{}{% Remove url except for @online
  \ifentrytype{unpublished}{}{% ou bien @unpublished
    \clearfield{url}
  }
  }
}

\renewbibmacro{in:}{%
    \ifentrytype{article}{}{\printtext{\bibstring{in}\intitlepunct}}}

\DeclareFieldFormat[article,periodical,inreference]{number}{\mkbibparens{#1}}
\DeclareFieldFormat[article,periodical,inreference]{volume}{\mkbibbold{#1}}
\renewbibmacro*{volume+number+eid}{%
    \printfield{volume}%
    \setunit*{\addthinspace}% NEW
    \printfield{number}%
    \setunit{\addcomma\space}%
    \printfield{eid}}

\DeclareFieldFormat[article,inbook,incollection]{title}{\enquote{#1}\addcomma} 

\usepackage{xspace}
\newcommand{\eme}{\up{\lowercase{e}}\xspace}

\newcommand{\numero}{n\textsuperscript o\kern+.2em}
\newcommand{\numeros}{n\textsuperscript {os}\kern+.2em}

  \DeclareMathSymbol{A}{\mathalpha}{operators}{`A}%
  \DeclareMathSymbol{B}{\mathalpha}{operators}{`B}%
  \DeclareMathSymbol{C}{\mathalpha}{operators}{`C}%
  \DeclareMathSymbol{D}{\mathalpha}{operators}{`D}%
  \DeclareMathSymbol{E}{\mathalpha}{operators}{`E}%
  \DeclareMathSymbol{F}{\mathalpha}{operators}{`F}%
  \DeclareMathSymbol{G}{\mathalpha}{operators}{`G}%
  \DeclareMathSymbol{H}{\mathalpha}{operators}{`H}%
  \DeclareMathSymbol{I}{\mathalpha}{operators}{`I}%
  \DeclareMathSymbol{J}{\mathalpha}{operators}{`J}%
  \DeclareMathSymbol{K}{\mathalpha}{operators}{`K}%
  \DeclareMathSymbol{L}{\mathalpha}{operators}{`L}%
  \DeclareMathSymbol{M}{\mathalpha}{operators}{`M}%
  \DeclareMathSymbol{N}{\mathalpha}{operators}{`N}%
  \DeclareMathSymbol{O}{\mathalpha}{operators}{`O}%
  \DeclareMathSymbol{P}{\mathalpha}{operators}{`P}%
  \DeclareMathSymbol{Q}{\mathalpha}{operators}{`Q}%
  \DeclareMathSymbol{R}{\mathalpha}{operators}{`R}%
  \DeclareMathSymbol{S}{\mathalpha}{operators}{`S}%
  \DeclareMathSymbol{T}{\mathalpha}{operators}{`T}%
  \DeclareMathSymbol{U}{\mathalpha}{operators}{`U}%
  \DeclareMathSymbol{V}{\mathalpha}{operators}{`V}%
  \DeclareMathSymbol{W}{\mathalpha}{operators}{`W}%
  \DeclareMathSymbol{X}{\mathalpha}{operators}{`X}%
  \DeclareMathSymbol{Y}{\mathalpha}{operators}{`Y}%
  \DeclareMathSymbol{Z}{\mathalpha}{operators}{`Z}%

\let\epsilon\varepsilon
\let\phi\varphi
\let\leq\leqslant
\let\geq\geqslant

\let\phi\varphi

\def\C{{\mathbf C}}
\def\R{{\mathbf R}}

\def\Z{{\mathbf Z}}

\title{Roger Godement et les fonctions de type positif}
\author{Alexandre Afgoustidis}
\address{CNRS \& Institut \'Elie Cartan de Lorraine, Nancy \& Metz, France}
\email{alexandre.afgoustidis@math.cnrs.fr}
\begin{document}
\selectlanguage{french}
\begin{abstract}
Ce texte, écrit pour la Gazette de la Société mathématique de France, évoque les fonctions de type positif et leur histoire avant 1950 ; on y  présente notamment des extraits de lettres écrites par Roger Godement, qui leur consacra sa thèse en 1946. 
\end{abstract}

% This is an expository paper on the early history of positive-definite functions, written for the ``Gazette de la Société mathématique de France''. It contains pictures of letters written by Roger Godement, whose 1946 thesis about positive-definite functions on groups was . 

\maketitle

\noindent

Les fonctions de type positif apparaissent en analyse et en probabilités vers 1910. Trente ans plus tard, elles contribuent à ouvrir la voie vers l'étude des représentations de dimension infinie des groupes topologiques localement compacts. C'est sur cela que porte la thèse de Godement, en 1946.  

On a retrouvé récemment des lettres que Godement envoyait entre 1945 et 1955 à Henri Cartan, qui fut son «\,directeur\,» de thèse. Le comité de rédaction de la \emph{Gazette} m'a proposé de présenter quelques extraits de ces lettres, en les replaçant dans leur contexte mathématique. C'est l'occasion de donner un aperçu de l'ambiance de cette époque, mais aussi d'évoquer la circulation des idées mathématiques dans les années 1930 et 1940 ; et de voir comment une notion apparue dans des questions d'analyse classique en vint à jouer un rôle important dans l'évolution de la théorie des représentations. 

\section{Les fonctions de type positif de 1909 à 1932}

\subsection{Définition et premier exemple}\label{sec:defs_generales}
Soient~$E$ un ensemble et~$K$ une fonction de~$E \times E$ dans~$\C$. Pour tout entier $n \geq 1$ et pour tout $n$-uplet $(x_1, \dots, x_n)$ d'éléments de~$E$, formons la matrice 
\begin{equation} \label{matrice_M} M_K(x_1, \dots, x_n)= (K(x_i,x_j))_{1 \leq i,j \leq n}.\end{equation} On dit que~$K$ est \emph{de type positif} si cette matrice est toujours hermitienne positive, quels que soient l'entier $n \geq 1$ et le $n$-uplet $(x_1, \dots, x_n)$.

Soient~$G$ un groupe et $f$ une fonction de~$G$ dans~$\C$. On dit que~$f$ est \emph{de type positif} si la fonction \mbox{$K_f\colon (x,y) \mapsto f(x^{-1}y)$}, de $G \times G$ dans~$\C$, est de type positif au sens ci-dessus. On notera alors $M_f(x_1, \dots, x_n)$ pour la matrice~\eqref{matrice_M}, plutôt que $M_{K_f}(x_1, \dots, x_n)$.

Par exemple, si~$f$ est une fonction de~$\R$ dans~$\C$, alors~$f$ est de type positif si pour tout $n \geq 1$ et pour tout $n$-uplet $(x_1, \dots, x_n)$ de nombres réels, la matrice 
\[ M_f(x_1, \dots, x_n)= (f(x_j-x_i))_{1 \leq i,j \leq n} \]
est hermitienne et positive. 

Donnons un premier exemple de fonction de type positif sur~$\R$. Fixons un nombre réel~$\lambda$ et considérons la fonction $e_{\lambda}\colon x \mapsto \exp(i\lambda x)$.   Pour tout $n \geq 1$ et tout \mbox{$(x_1, \dots, x_n) \in \R^n$}, la matrice $M_{e_\lambda}(x_1, \dots, x_n)$ est $(\overline{\exp(i\lambda x_k)} {\exp(i\lambda x_{\ell})})_{1 \leq k, \ell \leq n}$.  Pour voir qu'elle est hermitienne positive, il suffit de remarquer que c'est une \emph{matrice de Gram}. Rappelons en effet que si ${H}$ est un espace hilbertien complexe et si $(u_1, \dots, u_n)$ est une famille finie de vecteurs de~$E$, alors la matrice 
\[ \mathrm{Gram}_{{H}}(u_1, \dots, u_n) = (\langle u_i, u_j \rangle_{{H}})_{1 \leq i,j \leq n}\]
est hermitienne positive. Si~${H} = \C$ avec sa structure usuelle d'espace hilbertien, et si nous posons $u_j = \exp(-i \lambda x_j)$ pour tout $j$, alors  $\mathrm{Gram}_{{H}}(u_1, \dots, u_n)  = M_{e_\lambda}(x_1, \dots, x_n)$. Ainsi, la fonction~$e_\lambda$ est de type positif. 

%Toute combinaison linéaire à coefficients positifs de fonctions de type positif est encore une fonction de type positif, puisque les matrices hermitiennes positives $n \times n$ forment un cône convexe dans l'espace des matrices $n \times n$ à coefficients complexes. On en déduit, par exemple, que la fonction $\cos = \frac12(e_1 + e_{-1})$ est de type positif. 

\subsection{Mercer et les noyaux de type positif}\label{sec:mercer}

Le premier à évoquer les fonctions de type positif semble être James Mercer, dans un texte de~\cite*{Mercer} sur la théorie des équations intégrales. Pour en décrire le contexte, fixons des nombres réels $a$ et $b$ avec $a<b$ et notons $E$ l'espace des fonctions continues de $I=[a,b]$ dans~$\R$. Fixons aussi des fonctions continues \mbox{$K\colon I\times I \to \R$} et \mbox{$\varphi\colon I \to \R$}. Dans les premières années du \textsc{XX}\eme siècle, Ivar Fredholm~\parencite*{Fredholm} puis David Hilbert~\parencite*{Hilbert_1904} étudient les équations de la forme $f + \int_{I} K(\cdot,y) f(y) dy = \varphi$, où l'inconnue est une fonction $f \in E$. Hilbert fait un grand pas en reliant l'étude de cette équation à la «\,réduction\,» de la forme bilinéaire 
\begin{equation}\label{fq_mercer} (f, g) \mapsto \int_{I\times I} K(x,y) f(x) g(y) dxdy\end{equation}
sur $E \times E$. C'est le point de départ de la « théorie spectrale » en dimension infinie.

Alors que les premiers travaux de Hilbert sur ces questions sont tout récents, \textcite{Mercer} demande : parmi les noyaux $K\colon I\times I \to \R$ qui sont symétriques (c'est-à-dire vérifient $K(x,y) = K(y,x)$ pour tout $(x,y)$),  quels sont ceux pour lesquels la forme bilinéaire~\eqref{fq_mercer} est positive ? Il découvre que ce sont précisément les noyaux qui sont de type positif au sens du \S\,{}\ref{sec:defs_generales}. Il montre ensuite que pour ces noyaux, la forme~\eqref{fq_mercer} admet une bonne «\,réduction\,» au sens de Hilbert : dans la terminologie actuelle, cela revient à diagonaliser l'opérateur $f \mapsto \int_{I} K(\cdot, y) f(y) dy$ dans une base hilbertienne de~$L^2(I;\R)$. J'évoquerai au \S\,{}\ref{sec:stochastiques} la postérité des idées de Mercer dans l'étude des processus stochastiques. 

Pour des précisions sur les résultats de Mercer et de ses contemporains, et pour bien d'autres aspects de l'histoire des fonctions de type positif, on pourra consulter le texte de James~Stewart~\parencite*{Stewart}. 

%\small 
%
%
%\begin{itemize}
%\item pour montrer que si~\eqref{fq_mercer} est positive, alors la matrice $(K(x_i,x_j))_{1 \leq i,j \leq n}$ est positive quels que soient $n \geq 1$ et $(x_1, \dots, x_n) \in E^n$ , son argument revient à approcher la distribution $\sum \delta_{x_i}$ (somme des masses de Dirac en les $x_i$) par des fonctions affines par morceaux ; 
%\item pour montrer que si $K$ est de type positif, alors la forme quadratique \eqref{fq_mercer} est positive, il utilise les travaux de Fredholm et Hilbert qui reviennent à diagonaliser l'opérateur compact $f \mapsto \int_{I} K(\cdot, y) f(y) dy$.
%\end{itemize}  
%\normalsize

\subsection{Des exemples probabilistes}

Si Mercer est l'un des premiers à évoquer les fonctions de type positif (et semble créer le terme), la source la plus importante pour les développements en théorie des groupes est  à chercher plus près des probabilités que des équations intégrales. 

Si~$\mu$ est une mesure de probabilité sur~$\R$, convenons que sa \emph{transformée de Fourier} est la fonction~$\widehat{\mu}\colon \R \to \C$  définie par $\widehat{\mu}(\lambda) = \int_{\R} e^{-i\lambda x} d\mu(x)$. Dans le langage probabiliste, la fonction $\lambda \mapsto \hat{\mu}(- \lambda)$ est la fonction caractéristique d'une variable aléatoire distribuée selon la loi~$\mu$.

Si $(\lambda_1, \dots, \lambda_n)$ est un $n$-uplet de nombres réels avec $n \geq 1$, alors  les fonctions  \mbox{$e_{\lambda_j}\colon \R\to\C$} fournissent des éléments de l'espace hilbertien~${L}^2(\R;\mu)$, et on constate aussitôt que 
\[ M_{\hat{\mu}}(\lambda_1, \dots, \lambda_n) = \mathrm{Gram}_{{L}^2(\R;\mu)}(e_{\lambda_1}, \dots, e_{\lambda_n}). \]
La transformée de Fourier $\hat{\mu}$ est donc une fonction de type positif sur~$\R$. De plus, elle est continue (par convergence dominée) et vérifie $\hat{\mu}(0)=1$.

On obtient ainsi beaucoup de fonctions de type positif « naturelles » sur~$\R$ ; la fonction~$e_\lambda$ du \S\,{}\ref{sec:defs_generales} correspond au cas où $\mu$ est la masse de Dirac en~$\lambda$.  

La notion de fonction de type positif sur~$\R$ (en une variable), et le procédé ci-dessus pour en construire, semblent apparaître pour la première fois dans un texte de Maximilian Mathias~\parencite*{Mathias}. Je n'ai pas connaissance d'indications laissant penser que Mathias connaissait les travaux sur les équations intégrales du \S\,{}\ref{sec:mercer}. Il parle de «\,positiv-definite Funktion\,», ce qui est probablement à l'origine de la terminologie anglophone actuelle (qui me semble moins heureuse, à mon avis, que celle de~Mercer).

Mathias cherchait à adapter, pour  la transformation de Fourier dans~$\R$, des observations sur les séries de Fourier faites autour de 1910. 
Considérons $\mathbb{S}^1 = \R/\Z$, vu comme le cercle unité dans~$\C$ ; à toute mesure positive~$\mu$ sur~$\mathbb{S}^1$, on peut associer la suite $\hat{\mu}=(\hat{\mu}_n)_{n \in \Z}$ des coefficients de Fourier de~$\mu$, définie par $\hat{\mu}_n = \int_{\mathbb{S}^1} \overline{z}^n d\mu(z)$ pour $n \in \Z$. Quelles sont les suites de nombres complexes que l'on peut obtenir ainsi ? C'est le « problème des moments trigonométrique~», résolu par Gustav Herglotz~\parencite*{Herglotz}. 

\medskip

\begin{theo}\label{th:herglotz}
Soit~$(u_n)_{n \in \Z}$ une suite de nombres complexes. Les conditions suivantes sont équivalentes : 
\begin{enumerate}
\item[\textup{(i)}] Il existe une mesure positive de masse finie~$\mu$ sur~$\mathbb{S}^1$ vérifiant $u_n=\hat{\mu}_n$ pour tout $n \in \Z$~\textup{;}
\item[\textup{(ii)}] La suite $(u_n)_{n \in \Z}$ est de type positif\textup{.}
\end{enumerate}
\end{theo}

\medskip

La notion de suite de type positif venait d'apparaître, implicitement, chez Otto Toeplitz~\parencite*{Toeplitz}, pour reformuler algébriquement des conditions qu'obtenait Constantin Carathéodory sur les coefficients de Fourier de certaines fonctions analytiques d'une variable complexe. 

Mathias prouve que si~$f\colon \R \to \C$ est une fonction continue intégrable dont la transformée de Fourier est de type positif, alors~$f(x) \geq 0$ pour tout $x \in \R$  ; mais il n'obtient pas tout à fait l'analogue du théorème~\ref{th:herglotz} pour la transformation de Fourier sur~$\R$. Le mérite en revient à Salomon Bochner, dans ses leçons (\emph{Vorlesungen}) sur les intégrales de Fourier~\parencite*{Bochner_Vorlesungen} : 

\medskip

\begin{theo}\label{th:bochner}
Soit~$f$ une fonction de~$\R$ dans~$\C$. Les conditions suivantes sont équivalentes : 
\begin{enumerate}
\item Il existe une mesure de probabilité~$\mu$ sur~$\R$ vérifiant $f = \hat{\mu}$ \textup{;} 
\item La fonction~$f$ est continue, de type positif, et vérifie $f(0)=1$\textup{.} 
\end{enumerate}
\end{theo}

\section{Des probabilités aux groupes}

Vers 1932, les fonctions de type positif sont donc fermement ancrées dans l'analyse classique et, dirions-nous aujourd'hui, dans les probabilités. Comment cette notion migre-t-elle vers la théorie des représentations de groupes ?

Les représentations dont il est question ici sont des représentations \emph{unitaires}. Soit~$G$  un groupe topologique. Si~$H$ est un espace hilbertien complexe et si~$\mathcal{U}(H)$ désigne le groupe des automorphismes unitaires de~$H$, alors une  \emph{représentation unitaire de~$G$ dans~$H$} est un morphisme $\pi\colon G \to \mathcal{U}(H)$ tel que l'application $(g, x) \mapsto \pi(g)x$ soit continue de $G \times H$ dans~$H$.

\subsection{Le lien entre groupes, analyse harmonique et mécanique quantique}

Entre 1926 et~1932, des liens profonds apparaissent entre des thèmes auparavant éloignés\footnote{Le lien entre (a) et (c) était bien sûr connu ; mais avant 1926, il ne semble pas que les autres l'aient été.} : 
\begin{itemize}
\item[(a)] l'analyse harmonique (disons, de Fourier) ;
\item[(b)] la théorie des représentations de groupes ;
\item[(c)] la théorie spectrale et l'analyse fonctionnelle ;
\item[(d)] la Mécanique quantique.
\end{itemize}
L'histoire de ces rapprochements est fascinante. Je ne peux pas la retracer ici,  mais Bourbaki, par exemple, l'a fait récemment~\parencite*{TS345}. 

Disons simplement que c'est Hermann Weyl qui tisse les liens entre les représentations de groupes et les autres sujets ci-dessus. Il opère en 1926--1927 la jonction entre représentations des groupes compacts\footnote{L'étude des représentations des groupes compacts était aussi jeune que la Mécanique quantique : jusqu'en novembre 1924 (travaux de Schur, puis Weyl), on ne s'intéressait aux représentations que pour les groupes finis.}, analyse harmonique et théorie spectrale \parencite{Peter_Weyl}, guidé par le lien entre l'analyse de Fourier sur le cercle~$\mathbb{S}^1$ et les représentations du groupe commutatif compact correspondant\footnote{Si l'on voit~$\mathbb{S}^1$ comme le cercle unité dans~$\mathbb{C}$, alors toute représentation irréductible de~$\mathbb{S}^1$ sur un espace hilbertien~$H$ est de la forme $e^{i\theta} \mapsto e^{in\theta} \mathrm{Id}_H$, où~$n$ est un élément de~$\Z$. Les fonctions $e^{i\theta} \mapsto e^{in\theta}$ forment une base hilbertienne de ${L}^2(\mathbb{S}^1)$, et l'analyse harmonique s'intéresse beaucoup à cette base. Pour décrire sommairement ce que Weyl fait de cette idée dans le cas d'un groupe compact~$G$ arbitraire, mentionnons que si~$\mu$ est une mesure sur~$G$ invariante par translations, alors Peter et Weyl~\parencite*{Peter_Weyl}  construisent une base hilbertienne de $L^2(G; \mu)$ à partir des représentations unitaires irréductibles de~$G$ --- plus précisément, une base formée de coefficients matriciels de représentations irréductibles (ces coefficients matriciels sont définis ci-dessous, au \S\,\ref{sec:riesz_stone})}. Sa généralisation à tous les groupes compacts utilise de manière essentielle  la théorie spectrale telle qu'étudiée par Hilbert et ses élèves. Aussitôt après, dans des textes retentissants~(\cite{Weyl_Gruppen1}, \cite*{Weyl_Gruppen}),  il montre ce que les représentations peuvent apporter à la jeune Mécanique quantique. Voir les récits de Thomas Hawkins~\parencite*[ch.~4]{Hawkins} et d'Armand Borel~\parencite*{Borel}.
 
Quant au lien entre théorie spectrale et Mécanique quantique, on sait le rôle qu'eut John von~Neumann dans sa formulation, en définissant les espaces hilbertiens et les opérateurs partiels (habituellement dits « non bornés »), en dévoilant leur rôle pour la formulation de la théorie quantique, et en prouvant le théorème spectral pour les opérateurs partiels auto-adjoints. Rappelons sommairement que dans la formulation de von~Neumann, par exemple dans les \emph{Mathematische Grundlagen der Quantenmechanik}~\parencite*{vonNeumann_Grundlagen}, l'espace des états d'un système quantique est représenté par un espace hilbertien complexe~$H$, et les grandeurs physiques observables par des opérateurs partiels auto-adjoints sur~$H$ ; le théorème spectral donne accès aux spectres de ces opérateurs, et ainsi aux nombres réels représentant les résultats de mesures des grandeurs en question. 

Ces développements suscitent évidemment beaucoup de travaux mêlant les thèmes \mbox{(a)--(d)}. 

Un exemple frappant est le \emph{théorème de Stone} sur les représentations unitaires du groupe des translations de~$\R$. Weyl avait remarqué que si~$H$ est un espace hilbertien complexe et si~$u$ est un endomorphisme continu auto-adjoint de~$H$, alors l'application $t \mapsto \exp(itu)$, de~$\R$ dans $\mathcal{U}(H)$, est une représentation unitaire du groupe~$\R$ ; certains opérateurs venus de la mécanique quantique fournissaient donc des représentations de groupe. Il avait espéré que toute représentation unitaire de~$\R$ soit en fait de la forme $t \mapsto \exp(itu)$ où~$u$ n'est plus nécessairement un endomorphisme continu de~$H$, mais un \emph{opérateur partiel auto-adjoint} sur~$H$ (tel qu'étudié par von~Neumann dans le contexte ci-dessus)\footnote{On dit que~$u$ est le \emph{générateur infinitésimal} de la représentation donnée ; la théorie moderne des probabilités connaît bien, elle aussi, les générateurs infinitésimaux de (semi-)groupes à un paramètre d'opérateurs stochastiques. }.  Presque aussitôt, c'est Marshall Stone~\parencite*{Stone_Weyl}, en marge de la rédaction de son grand livre sur la théorie spectrale hilbertienne\footnote{Le livre de Stone~\parencite*{Stone_Livre} fut la première référence «\,canonique\,» pour le point de vue post-hilbertien sur la théorie spectrale, à la suite des travaux de von~Neumann et de Stone lui-même. }, qui démontre la «\,conjecture\,» de Weyl à l'aide du théorème spectral. 

Le théorème de Stone est l'un des premiers résultats significatifs sur les représentations unitaires d'un groupe non compact ; comme on le voit, il est motivé par la mécanique quantique, et en lien étroit avec la théorie spectrale hilbertienne.

Les liens entre nos thèmes~(a) à~(d) étaient donc très porteurs à l'époque où Bochner prépare ses \emph{Vorlesungen} et démontre le théorème~\ref{th:bochner}. Il a de bonnes raisons de se tenir au courant des liens qu'on vient d'évoquer : Bochner était l'un des meilleurs spécialistes des «\,fonctions presque périodiques\,» d'Harald Bohr, qui ont connu une grande vogue à peu près à la même période, et qui apparaissent souvent dans les textes de Weyl et von~Neumann\footnote{Soient~$G$ un groupe et~$f\colon G \to \C$ une fonction continue bornée. On dit que~$f$ est \emph{presque périodique} si, dans l'espace des fonctions continues bornées de~$G$ dans~$\C$, les ensembles de translatées $\{(g \mapsto f(ag)), a \in G\}$ et $\{(g \mapsto f(ga)), a \in G\}$, sont relativement compacts pour la topologie de la convergence uniforme. Cette notion est due à von~Neumann \parencite*{vonNeumann_almostperiodic} et généralise la notion de fonction presque périodique sur~$\R$ d'Harald Bohr~\parencite*{Bohr}. Comme de très nombreux travaux des années 1925--1933 avaient montré que les fonctions presque périodiques sur~$\R$ ont des propriétés remarquables du point de vue de l'analyse harmonique, la généralisation à un groupe quelconque a soulevé de grands espoirs. Mais les fonctions presque périodiques sont aujourd'hui presque oubliées. Il n'est probablement pas absurde de dire, comme  Weil et Chevalley dans leur nécrologie mathématique de Weyl~\parencite{Chevalley_Weil}, que «\,leur rôle principal a été de préparer le point de vue moderne sur les groupes localement compacts\,».}. 

En~1933, dans le noir contexte que l'on sait, Bochner et Weyl quittent l'Europe et acceptent des postes à l'IAS de Princeton. Ils y rejoignent von~Neumann, qu'ils semblent déjà bien connaître tous les deux. Il est probable que Bochner et von~Neumann aient beaucoup discuté des liens entre analyse harmonique, théorie des groupes et analyse fonctionnelle. Ils étudieront ensemble les fonctions presque périodiques sur un groupe quelconque~\parencite*{Bochner_vonNeumann}.

\subsection{Riesz, Bochner et le théorème de Stone}\label{sec:riesz_stone}

C'est aussi en~1933  qu'apparaît le premier lien entre fonctions de type positif et représentations de groupes. Il est établi par Frédéric Riesz, et indépendamment par Bochner lui-même, à propos du théorème de Stone \parencite{Riesz_1933, Bochner_1933}. C'est l'observation suivante. 

Soit~$G$ un groupe topologique. Considérons une représentation unitaire $\pi\colon G \to \mathcal{U}({H})$ de~$G$  dans un espace hilbertien complexe~${H}$.

À tout vecteur~$\xi$ de~${H}$, associons une fonction $c_{\pi, \xi}\colon G \to \C$ en posant $c_{\pi, \xi}(g)= \langle \xi, \pi(g)\xi \rangle_{{H}}$. On dit que $c_{\pi, \xi}$ est le \emph{coefficient matriciel diagonal}\footnote{Les coefficients non nécessairement diagonaux sont les fonctions de la forme $g \mapsto  \langle \xi, \pi(g)\eta \rangle_{{H}}$ pour $\xi, \eta \in H$.} de~$\pi$ associé au vecteur~$\xi$. Cette fonction est toujours continue et de type positif : en effet, si $(g_1, \dots, g_n)$ est une famille finie de points de~$G$, alors  $M_{c_{\pi,\xi}}(g_1, \dots, g_n)=\mathrm{Gram}_H(\pi(g_1)\xi, \dots, \pi(g_n)\xi)$. 

Riesz et Bochner observent ce qui précède dans le cas où~$G$ est le groupe des translations de~$\R$. Et ils montrent que cela permet de \emph{déduire le théorème de Stone de celui de Bochner}. 

Voici une idée vague de l'argument de Riesz. Pour tout $\xi \in H$, on vient de voir que $c_{\pi,\xi}$ est une fonction de type positif sur~$\R$. D'après le théorème de Bochner, c'est la transformée de Fourier d'une mesure de probabilité~$\mu_{\xi}$ sur~$\R$. Reste à déduire des mesures $\mu_\xi$, $\xi \in H$, un opérateur partiel auto-adjoint $u$ sur~$H$, et à montrer qu'on~a $\pi(u)=\exp(itu)$ pour tout $t \in \R$. Cela se fait par le truchement du théorème spectral : dans l'une de ses versions classiques, ce dernier prend la forme d'une bijection entre  opérateurs partiels auto-adjoints sur~$H$ et  \emph{résolutions de l'identité} de~$H$ (sortes de mesure de probabilité sur~$\R$ à valeurs dans les orthoprojecteurs de~$H$). Riesz construit une telle résolution à partir des $\mu_\xi$, et en déduit le théorème de Stone. %Voir l'exposé de Cartan~\parencite*{Julia} au séminaire Julia. 

\subsection{Le séminaire Julia de 1934--1935}\label{sec:julia}

Le séminaire Julia, organisé entre 1933 et 1939, fut l'un des premiers séminaires de mathématiques en France, après le séminaire Hadamard organisé entre 1920 et 1937. Michèle Audin en a étudié et publié les archives~\parencite*{Julia_1}. En forçant beaucoup le trait\footnote{Michèle Audin explique en détail pourquoi cette vague parenté s'applique au moins aussi bien au séminaire Hadamard \parencite*[\S\,4.3]{Julia_1}. }, on peut y penser comme à une sorte d'ancêtre du séminaire Bourbaki : il se tenait à l'IHP, on y donnait des exposés sur l'actualité mathématique, et ces exposés étaient rédigés. C'est d'ailleurs en marge de ce séminaire que s'organisent les premières discussions sur un projet de «\,traité d'analyse\,», futur \emph{Éléments de mathématique}. 

Pour l'année 1934--1935, le thème choisi est : «\,espace de Hilbert\,»~\parencite*{Julia_2}. Les orateurs des trois dernières séances sont André Weil, John von~Neumann et Henri Cartan. Weil parle le 8 avril des fonctions presque périodiques. L'exposé du 6 mai est donné par von~Neumann ; annoncé à la séance précédente comme portant sur les représentations de groupes, il porte finalement sur les «\,anneaux d'opérateurs\,», avec  l'intention explicite d'appliquer les seconds aux premières. Le~20~mai, Cartan présente le mémoire de Riesz sur le théorème de Stone. 

Plusieurs des exposés donnés cette année-là font référence à un «\,mémorial (à paraître)\,» de Weil. Ce sera \emph{l'intégration dans les groupes topologiques et leurs applications}~\parencite*{Weil_Livre}, que Weil rédige en~1935 pour faire le point sur les résultats connus sur les représentations des groupes compacts, et sur les groupes localement compacts commutatifs et leurs caractères. Il y introduit les fonctions de type positif sur un groupe quelconque, et démontre l'analogue du théorème de Bochner sur un groupe abélien localement compact. Ce livre est toujours, à mon avis, une belle introduction au sujet. Weil évoquera plus tard~\parencite*{Weil_Oe} l'objectif implicite de l'ouvrage : 
\medskip

\begin{center}\begin{minipage}{0.8\textwidth}\small {\emph{Je l'avais entrepris surtout avec l'espoir d'ouvrir la voie à une généralisation de la théorie des représentations des groupes finis et des groupes compacts ; non seulement je n'atteignis pas la terre promise des représentations de dimension infinie, mais je m'arrêtai avant même de l'entrevoir. Je me décourageai trop tôt quand je vis que les coefficients des représentations de degré fini des groupes de Lie simples non compacts ne sont pas de carré intégrable [...]. % ; cela empêche, si on se propose d'explorer les représentations dans $L^2$, d'avancer du connu vers l'inconnu. 
Peut-être, pour une démarche aussi hardie, fallait-il l'esprit sans préjugé d'un physicien ; ce sont bien en effet Dirac, Wigner, Bargmann qui ont ouvert la voie à cet égard, à propos du groupe de Lorentz.} }\end{minipage}\end{center}
\normalsize
\medskip
Les «\,physiciens\,» en question sont motivés par l'idée suivante, évoquée par Dirac dès 1936 : si~$H$ est un espace hilbertien modélisant les états d'un système quantique, et si le système se plie aux lois de la relativité restreinte qui sont invariantes par les groupes de Lorentz $\mathrm{SO}(3,1)$ et de Poincaré $\mathrm{SO}(3,1) \ltimes \R^4$, alors $H$ doit porter une représentation de l'un de ces groupes.  En~1939, sur une suggestion de von Neumann,  Wigner classifie les représentations irréductibles de $\mathrm{SO}(3,1) \ltimes \R^4$ ; c'est la première étude détaillée des représentations d'un groupe qui n'est ni compact ni abélien~\parencite*{Wigner}.

\section{La thèse et les lettres de~Godement}

Roger Godement, né en~1921, entre à l'École Normale en 1940. Il fait partie de la première promotion appelée à suivre l'enseignement d'Henri Cartan, qui vient d'être nommé à Paris. Six ans plus tard, en juillet 1946, Godement est le premier mathématicien à soutenir sa thèse sous la «\,direction\,» de Cartan (qui tenait aux guillemets).  Le titre est : \emph{les fonctions de type positif et la théorie des groupes}~\parencite*{Godement_these}. J'espère que le \S\,{}\ref{sec:julia} montre que ce choix de sujet n'est pas vraiment une surprise ; nous verrons comment le travail de Godement contribue à ouvrir le chemin qu'espérait Weil.

Avant la Seconde Guerre mondiale, Cartan vivait à Strasbourg et avait dû quitter précipitamment la ville lors de l'évacuation de septembre 1939. Il y revient en 1945 pour un détachement de deux ans. Godement reste à Paris et lui écrit de très nombreuses lettres, soigneusement conservées par Cartan. L'ouverture du fonds Henri Cartan conservé à l'Académie des sciences a permis à Christophe Eckes de redécouvrir ces lettres en octobre 2021. Malheureusement, nous n'avons que très peu de copies des lettres de Cartan à Godement.
 
\begin{figure}[h]
\begin{center}\includegraphics[width=0.6\linewidth]{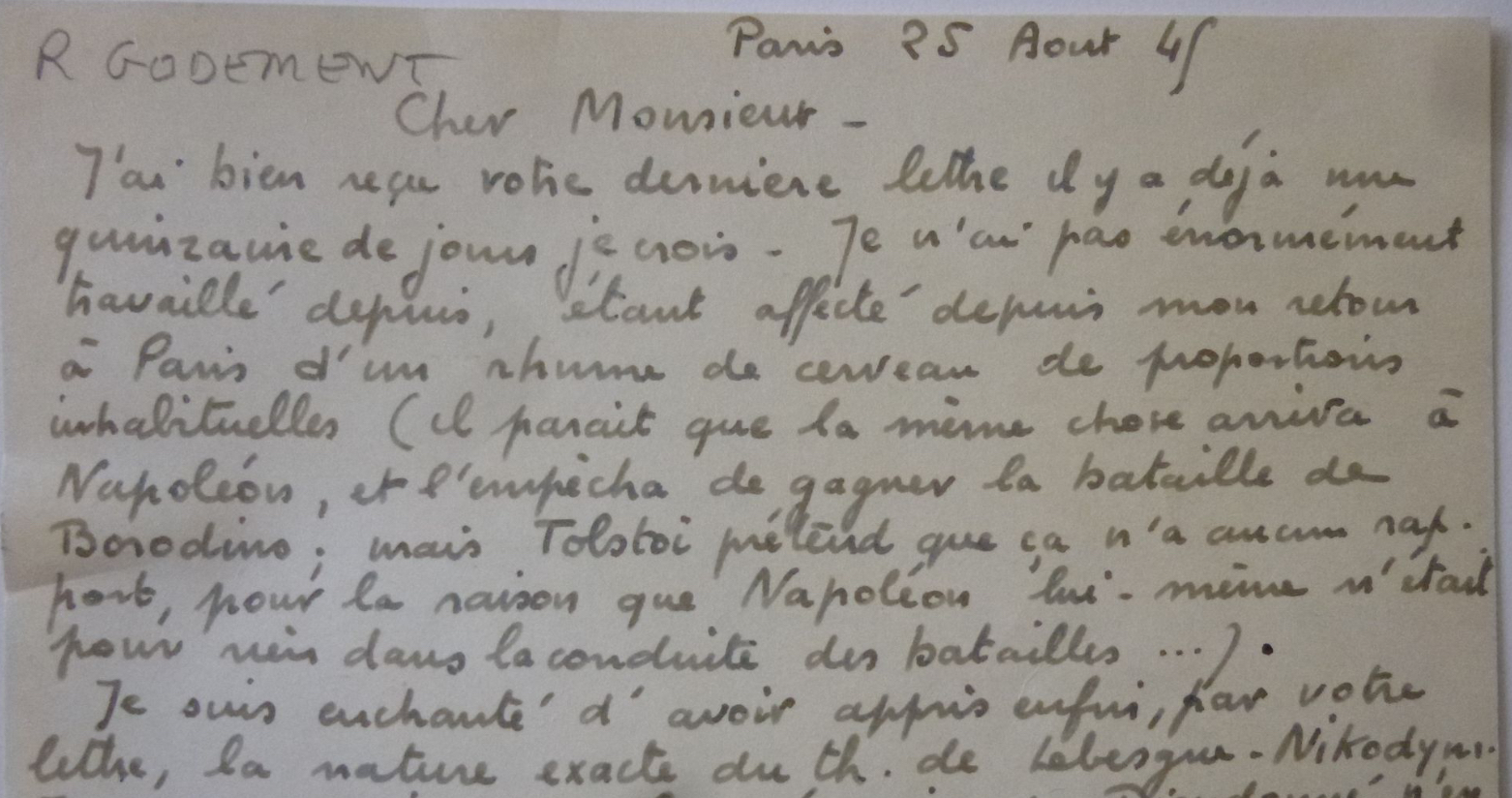}\end{center}
\caption{Début de la lettre du 25 août 1945}
\end{figure}
On dispose ainsi de 54 lettres écrites entre 1945 et 1955, dont 34 pour la seule année universitaire 1945--1946 qui s'achève sur la soutenance. 
Dans les lettres de 1945--1946, Godement parle surtout de mathématiques et de tracas matériels ou administratifs liés à sa fin de thèse.  La suite de la correspondance mêle discussions mathématiques, échanges liés à Bourbaki, et nouvelles de la famille ou des  postes universitaires.

Les extraits reproduits dans les \S\,\ref{sec:gelfand} à \ref{sec:circu} en sont un petit échantillon, choisi pour ce qu'il dit du travail mathématique de Godement ou de la circulation des idées à cette époque. Avant d'évoquer les lettres elles-mêmes, j'énonce deux résultats de la thèse de Godement (\S\,{}\ref{sec:real}--\ref{sec:partitions}).

\subsection{Le lien entre fonctions de type positif et  représentations unitaires}\label{sec:real}

Soit~$G$ un groupe localement compact. Nous avons vu au \S\,{}\ref{sec:riesz_stone}, suivant Riesz et Bochner (1933), que si \mbox{$\pi\colon G \to \mathcal{U}(H)$} est une représentation unitaire de~$G$ dans un espace hilbertien $H$, alors pour tout $\xi \in H$, la fonction $c_{\pi, \xi}\colon G \to \C$ définie par $c_{\pi, \xi}(g)=\langle \xi, \pi(g)\xi \rangle_{{H}}$ est continue et de type positif. Godement observe que, réciproquement, \emph{toutes les fonctions continues de type positif s'obtiennent à partir de représentations unitaires.}

Voyons comment Godement associe, à toute fonction continue  $f \colon G\to \C$ de type positif, une représentation unitaire de~$G$.  

Pour $x \in G$, notons $f_x$ la fonction $g \mapsto f(x^{-1} g)$. Soit $V_{f}$ l'ensemble des combinaisons linéaires (finies) de fonctions $f_x$, $x \in G$. 
Si $\varphi = \sum_{i \in I} \alpha_i f_{x_i}$ et $\psi = \sum_{j \in J} \beta_j f_{y_j}$ sont des éléments de~$V_{f}$, posons $\langle \varphi, \psi \rangle = \sum_{(i,j) \in I \times J} \alpha_i \overline{\beta_i} f(x_i^{-1} y_j)$. Dire que $f$ est de type positif, c'est précisément dire que cette formule définit\footnote{On a $\langle \varphi, \psi \rangle = \sum_{j \in J} \overline{\beta}_j \varphi(y_j)$, donc $\langle \varphi, \psi \rangle$ ne dépend que de $\varphi$ et pas de son écriture sous la forme $\sum_{i \in I} \alpha_i f_{x_i}$. Par ailleurs, si $f$ est de type positif, alors $f(g^{-1}) = \overline{f(g)}$ pour tout $g \in G$ puisque la matrice $\big(\begin{smallmatrix} f(1) & f(g) \\ f(g^{-1}) & f(1)\end{smallmatrix}\big)$ est hermitienne ; il en résulte que $\langle \phi, \psi \rangle = \overline{\langle \psi, \phi\rangle}$, ce qui prouve que $\langle \phi, \psi \rangle$ ne dépend aussi que de $\psi$. L'application $(\varphi, \psi) \mapsto  \langle \varphi, \psi \rangle$ ainsi obtenue est clairement bilinéaire et définit une forme hermitienne sur~$V_f$ d'après ce qui précède. Enfin, si $\varphi = \sum_{i \in I} \alpha_i f_{x_i}$ alors $\langle \varphi, \varphi \rangle = \sum \limits_{i \in I} \lvert \alpha_i \rvert^2$, donc la forme $\langle \cdot, \cdot \rangle$ est positive.} une \emph{forme hermitienne positive} sur~$V_f$. Cette forme est alors définie positive. En effet, supposons $\langle \varphi, \varphi \rangle = 0$ ; alors $\langle \varphi, \psi \rangle = 0$ pour tout $\psi$ d'après l'inégalité de Cauchy--Schwarz. Or, la définition de $\langle \cdot, \cdot \rangle$ implique aussitôt qu'on a $\langle \varphi, f_x \rangle = \varphi(x)$ pour tout $x \in E$, si bien que $\varphi=0$.

Soit~$H_f$ le complété de $V_f$ pour le produit scalaire~$\langle \cdot, \cdot \rangle$. Pour tout $g \in G$, soit $\pi(g)\colon V_f \to V_f$ l'application $\sum \alpha_i f_{x_i} \to \sum \alpha_i f_{x_i g}$ ; alors $\pi(g)$ induit un endomorphisme unitaire de~$H_f$. On constate alors que $\pi\colon G \to \mathcal{U}(H_f)$ est une représentation unitaire de~$G$ et vérifie $c_{\pi, f} = f$.   Godement a donc démontré : 

\medskip

\begin{theo}\label{th:realisation_hilbertienne}
Si $f$ est une fonction de~$G$ dans~$\C$, alors les conditions suivantes sont équivalentes : 
\begin{enumerate}
\item[\textup{(i)}] $f$ est continue et de type positif ;
\item[\textup{(ii)}] il existe une représentation unitaire $\pi$ de~$G$ dans un espace hilbertien~$H$, et il existe $\xi \in H$, tels que l'on ait $f=c_{\pi,\xi}$.  
\end{enumerate}
\end{theo}

\medskip

\normalsize 
Si l'on se donne une fonction~$f$ vérifiant~(i), alors la construction ci-dessus fournit un couple $(\pi, \xi)$ concret vérifiant $f = c_{\pi, \xi}$ et tel que le vecteur $\xi$ de~$H$ soit \emph{cyclique}, c'est-à-dire que $\pi(G) \xi$ engendre un sous-espace dense de~$H$. Godement remarque que si l'on se donne une autre représentation unitaire $\pi' \colon G \to \mathcal{U}(H')$ et un vecteur cyclique $\xi' \in H'$, alors on a $c_{\pi, \xi} = c_{\pi',\xi'}$ si et seulement s'il existe un isomorphisme $T\colon H \to H'$ qui envoie $\xi$ sur $\xi'$ et \emph{entrelace} $\pi$ et~$\pi'$, c'est-à-dire vérifie $T(\xi)=\xi'$ et $T \circ \pi(g) = \pi(g) T'$ pour tout $g \in G$. En particulier, dans le th.~\ref{th:realisation_hilbertienne}, la classe d'équivalence de la représentation~$\pi$ est uniquement déterminée par~$f$. 

\subsection{Partitions des fonctions de type positif}\label{sec:partitions}

Soient~$G$ un groupe localement compact et~$f$ une fonction continue de type positif sur~$G$. Supposons que~$f$ corresponde par le théorème~\ref{th:realisation_hilbertienne} à une représentation~$\pi$ de~$G$, uniquement déterminée à équivalence près. Peut-on donner une condition nécessaire et suffisante sur~$f$ pour que~$\pi$ soit irréductible ? Godement répond à cette question en introduisant la notion de \emph{partition} des fonctions de type positif.

Appelons \emph{partition de~$f$} toute famille finie $(f_i)_{i \in I}$ de fonctions de type positif vérifiant $f = \sum_{i \in I} f_i$. Disons que $f$ est \emph{élémentaire} si les seules partitions de~$f$ sont les partitions \emph{triviales} --- celles où tous les $f_i$ sont de la forme $\alpha_i f$, pour une famille $(\alpha_i)_{i \in I}$ de nombres réels positifs de somme~$1$. 

Supposons qu'on ait $f = c_{\pi, \xi}$, où $\pi\colon G \to \mathcal{U}(H)$ est une représentation unitaire et où \mbox{$\xi \in H$}. Supposons que~$\pi$ ne soit pas irréductible, c'est-à-dire qu'il existe un sous-espace fermé $H_1$ de~$H$ qui soit stable par tous les~$\pi(g)$ et différent de~$\{0\}$ et de~$H$. Soient~$\xi_1$ et $\xi_2$ les projections orthogonales de $\xi$ sur $H_1$ et l'orthogonal de~$H_1$, respectivement. Alors il résulte aussitôt des définitions qu'on a $f = c_{\pi, \xi_1} + c_{\pi, \xi_2}$. Godement vérifie que la partition $(c_{\pi, \xi_1}, c_{\pi, \xi_2})$ de~$f$ est non triviale. La fonction~$f$ n'est donc pas élémentaire.  Par contraposée, pour que $f$ soit élémentaire, il est nécessaire que $\pi$ soit irréductible. Godement démontre que c'est suffisant : 

\medskip

\begin{theo}\label{th:elementaires} Soit $f\colon G \to \C$ une fonction de type positif, et soit~$\pi$ une représentation unitaire de~$G$ dont $f$ est un coefficient matriciel diagonal ; alors $\pi$ est irréductible si et seulement si $f$ est élémentaire. 
\end{theo}

\subsection{La surprise de février 1946}\label{sec:gelfand}

Les théorèmes~\ref{th:realisation_hilbertienne} et~\ref{th:elementaires} sont aujourd'hui des résultats de base sur les représentations des groupes localement compacts. Godement les avait énoncés dans des notes aux Comptes-Rendus parues entre juillet 1945 et janvier 1946. Godement avait publié plusieurs autres notes à cette date : elles concernaient le cas des groupes abéliens, le lien entre fonctions de type positif et fonctions presque périodiques, et une généralisation du théorème ergodique de von~Neumann. 

Rebondissement le 6 février 1946 : Godement apprend «\,par hasard\,» que les théorèmes~\ref{th:realisation_hilbertienne} et~\ref{th:elementaires} ont été démontrés trois ans auparavant par Israel Gelfand et Dmitrii Raikov~\parencite*{Gelfand_Raikov}. 

\begin{figure}[h]
\begin{center}\includegraphics[width=0.6\linewidth]{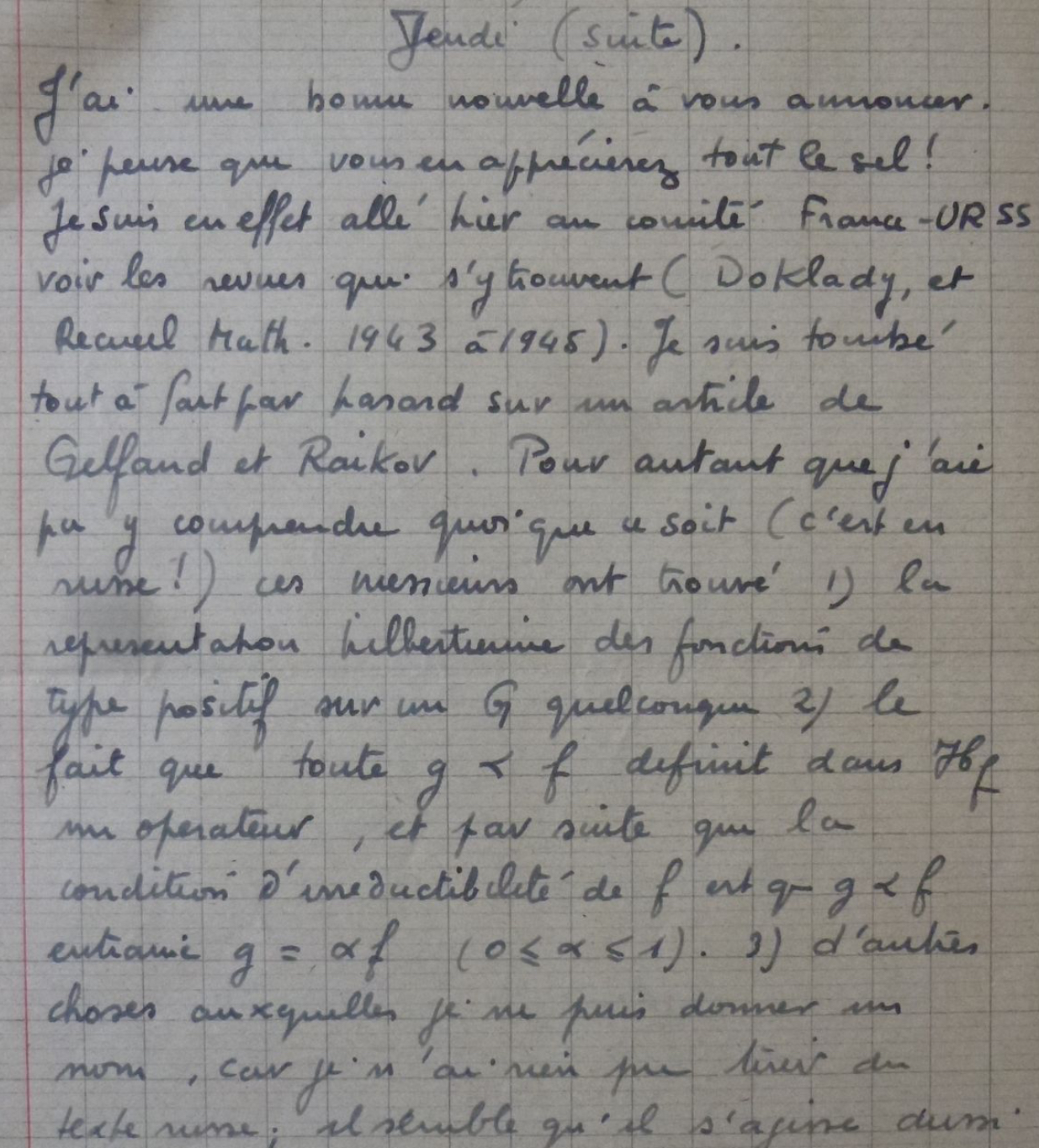}\end{center}
\caption{Début de la lettre du 7 février 1946}
\end{figure}

\subsection{Gelfand, Raikov et Krein}

En fait, Gelfand a commencé à explorer systématiquement, depuis 1940, les représentations des groupes localement compacts généraux. 
Comment a-t-il conçu ce projet ? Il n'avait apparemment pas connaissance du livre de Weil («\,\emph{il semble que mon livre n'atteignit l'Union Soviétique que lorsqu'Élie Cartan l'apporta à Moscou, ainsi que les premiers volumes de Bourbaki, à l'été 1945}\,»). 

Gelfand avait déjà à son actif, à la fin des années 1930, d'avoir jeté les bases de la théorie des algèbres d'opérateurs, et il avait éclairé par là certains aspects de l'analyse harmonique et des fonctions presque périodiques sur~$\R$.\footnote{La preuve que Gelfand et Raikov donnent pour le  th.~\ref{th:realisation_hilbertienne}, légèrement différente de celle du \S\,\ref{sec:real}, donnera d'ailleurs ce que l'on appelle aujourd'hui la « construction de Gelfand--Naimark--Segal » pour les $C^\ast$-algèbres. 
 } 

Il semble que l'idée de se tourner, pour les groupes non abéliens, vers les fonctions de type positif, puisse être reliée (\emph{via} \cite{Raikov}) à des travaux de Mark Krein. Ce dernier publie en~1940 une note reliant algèbres d'opérateurs et fonctions de type positif sur les groupes~\parencite*{Krein_note}. Krein connaissait bien les travaux de Mercer et le problème des moments, qu'il avait évoqués à plusieurs reprises dans les années 1930. C'était aussi un spécialiste de la convexité dans les espaces de Banach (sujet alors récent). Le célèbre théorème de Krein--Milman paraît d'ailleurs aussi en~1940 \parencite{Krein_Milman}. Les deux sujets se rejoignent si l'on observe que pour tout groupe localement compact~$G$, l'ensemble~$\mathrm{Pos}(G)$ des fonctions continues de type positif sur~$G$ est une partie \emph{convexe} de~$L^\infty(G)$.

Dans leur texte de 1943, Gelfand et Raikov définissent  les fonctions de type positif \emph{élémentaires} du \S\,\ref{sec:partitions}, et remarquent que ce sont les \emph{points extrémaux} de~$\mathrm{Pos}(G)$ --- ceux qui ne sont à l'intérieur d'aucun segment de droite contenu dans~$\mathrm{Pos}(G)$. 
Une application ingénieuse du théorème de Krein--Milman leur permet alors de montrer que les représentations unitaires irréductibles \emph{séparent les points} de~$G$ (ce résultat est aujourd'hui connu sous le nom de « théorème de Gelfand--Raikov ») : 
\medskip

\begin{theo}\label{th:gelfand_raikov} 
Si $g$ et $g'$ sont deux points distincts de $G$, alors il existe une représentation unitaire irréductible~$\pi$ de~$G$ telle qu'on ait $\pi(g)\neq \pi(g')$.
\end{theo}
\medskip
Découvrant le  «\,caractère rupinant des résultats\,» de Gelfand--Raikov, Godement s'aperçoit aussi qu'une bonne partie de ses propres idées n'est pas contenue dans leur mémoire, mais peut être précisée et simplifiée en adoptant leurs méthodes.\footnote{Disons un mot sur l'un des liens les plus simples entre fonctions de type positif et algèbres d'opérateurs, qui est à la base de ces simplifications. Soient~$G$  un groupe localement compact et~$\mu$ une mesure de Radon sur~$G$ invariante par les translations $h \mapsto gh$ pour $g \in G$ (une telle mesure « de Haar » existe toujours). Si~$\varphi$ et~$\psi$ sont des fonctions intégrables sur~$G$, définissons $\varphi \ast \psi$ comme la fonction $g \mapsto \int_{G} \varphi(h) \psi(gh^{-1})d\mu(h)$ sur~$G$, et  $\varphi^\dagger$ comme la fonction $g \mapsto \overline{\varphi(g^{-1})}$.  Ces formules induisent, sur l'espace de Banach~${L}^1(G)$, une loi de composition~$\ast$ et une involution $\dagger$. On dit qu'un élément~$\psi$ de~${L}^1(G)$ est \emph{positif} s'il peut s'écrire sous la forme $\varphi^\dagger \ast \varphi$, où~$\varphi \in L^1(G)$.  Si  $f\colon G \to \C$ est une fonction bornée, elle définit une forme linéaire $T_f$ sur~${L}^1(G)$ par la formule $T_f(\varphi)=\int_G f(g) \varphi(g) d\mu(g)$ ; la fonction~$f$ est alors \emph{de type positif} si et seulement si la forme linéaire~$T_f$ est \emph{positive sur tout élément positif} de~$L^1(G)$. Cela éclaire la définition des fonctions de type positif, et ouvre la voie à diverses généralisations --- très utiles pour la théorie des représentations et au-delà. }

\begin{figure}[h]
\begin{center} \includegraphics[width=0.5\linewidth]{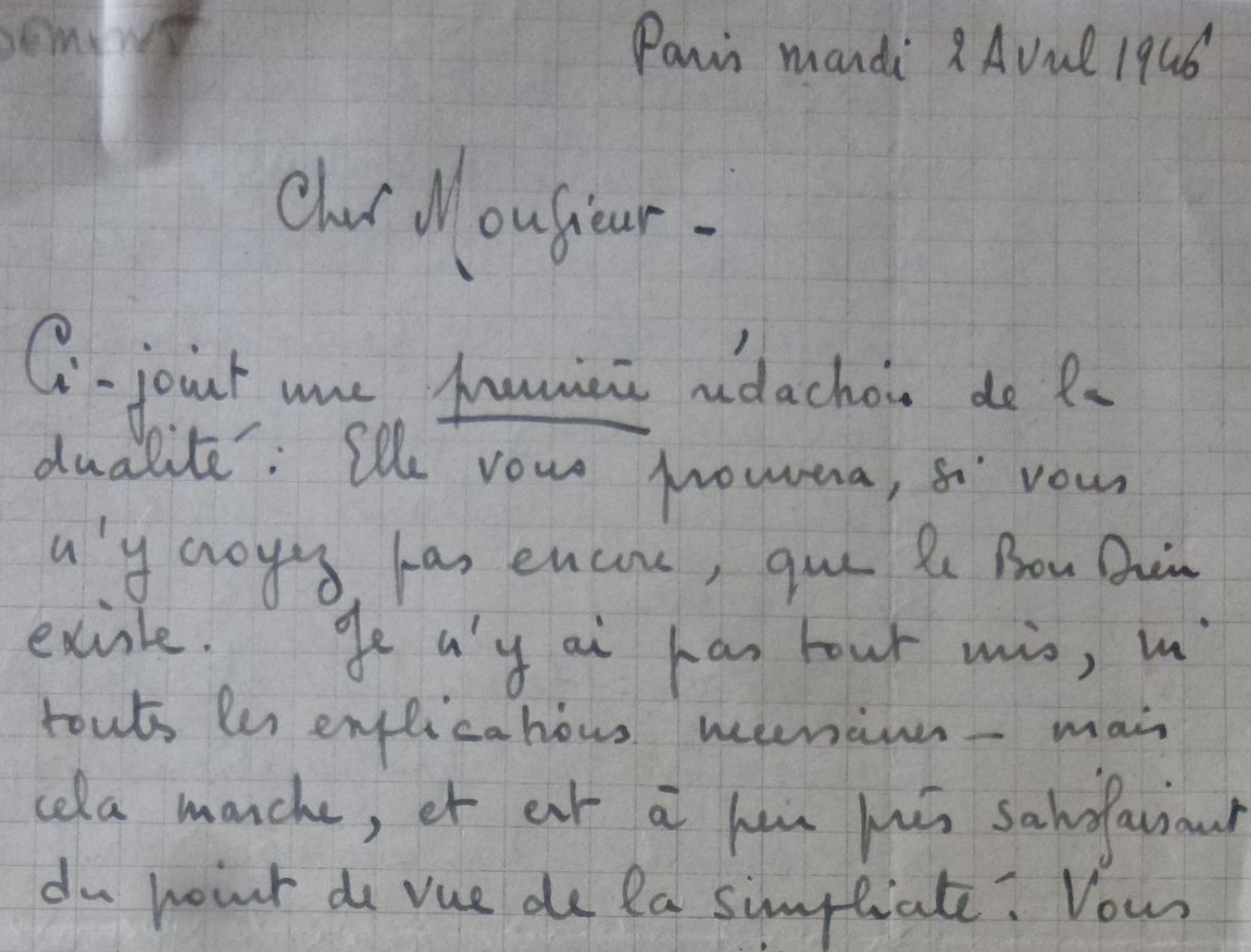}\end{center} 
\caption{Godement rassuré (avril 1946)}
\end{figure}

Moins de six mois plus tard,  le théorème de Krein--Milman devient l'un des outils de base de la version finale de sa thèse. Il utilise par exemple la convexité pour prouver que toute fonction continue de type positif sur~$G$ peut être approchée, uniformément sur tout compact, par des combinaisons linéaires de fonctions élémentaires (combinaisons qu'il appelle des «\,polynômes trigonométriques\,»).  De même, Cartan et Godement déduisent des méthodes de Gelfand et Raikov un exposé complet, nettement simplifié, de la dualité et de l'analyse harmonique pour les groupes abéliens~\parencite{Cartan_Godement}.

\begin{figure}[h]
\begin{center}\includegraphics[width=0.6\linewidth]{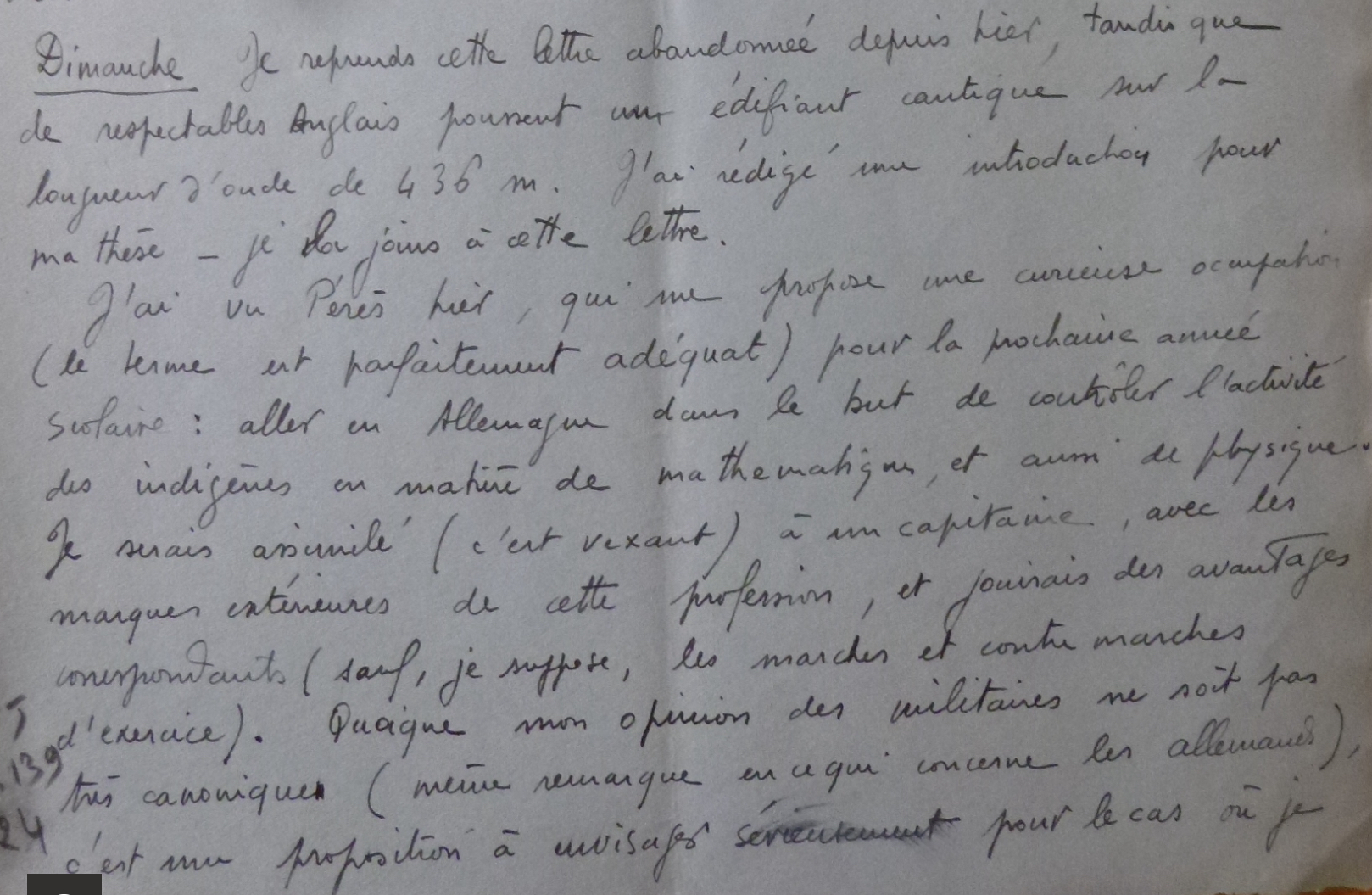}\end{center}
\caption{Un post-doc militaire ? (Mai 1946)}
\end{figure}
Tout est bien qui finit bien. Godement soutient sa thèse en juillet 1946 ; les membres de son jury sont Élie Cartan, Henri Cartan et Jean Favard. Il est aussitôt nommé à Nancy (ce qui semble mettre fin à des inquiétudes sur son avenir matériel, régulièrement exprimées en 1945--1946).

\subsection{Sur la circulation des idées}\label{sec:circu}

\begin{figure}
\begin{center}\includegraphics[width=0.55\linewidth]{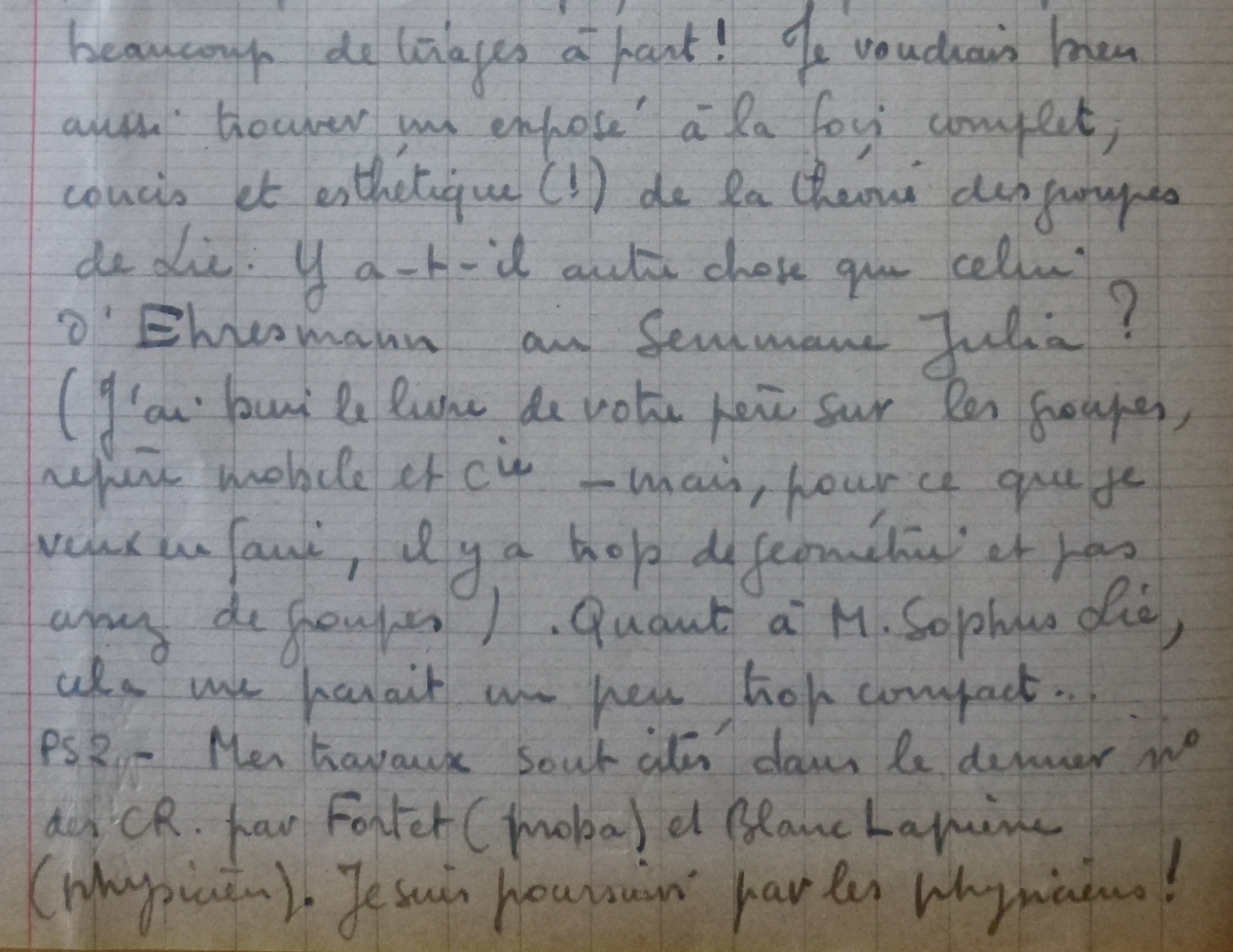}\end{center}
\caption{Poursuivi par les physiciens ! (Mars 1946)}\label{physiciens}
\end{figure}

D'autres extraits de la correspondance me semblent instructifs, à propos de l'évolution de la théorie des groupes à cette époque. Voyez par exemple cette lettre de mars 1946 où Godement cherche une référence «\,{}moderne\,» sur la théorie des groupes de Lie (Fig.~\ref{physiciens}).

Plus tard dans l'année paraîtra \emph{Theory of Lie groups}~\parencite*{Chevalley_Liegps}, de Chevalley --- qui  avait eu un rôle important dans la programmation du séminaire Julia.

\begin{figure}
\begin{center}\includegraphics[width=0.7\linewidth]{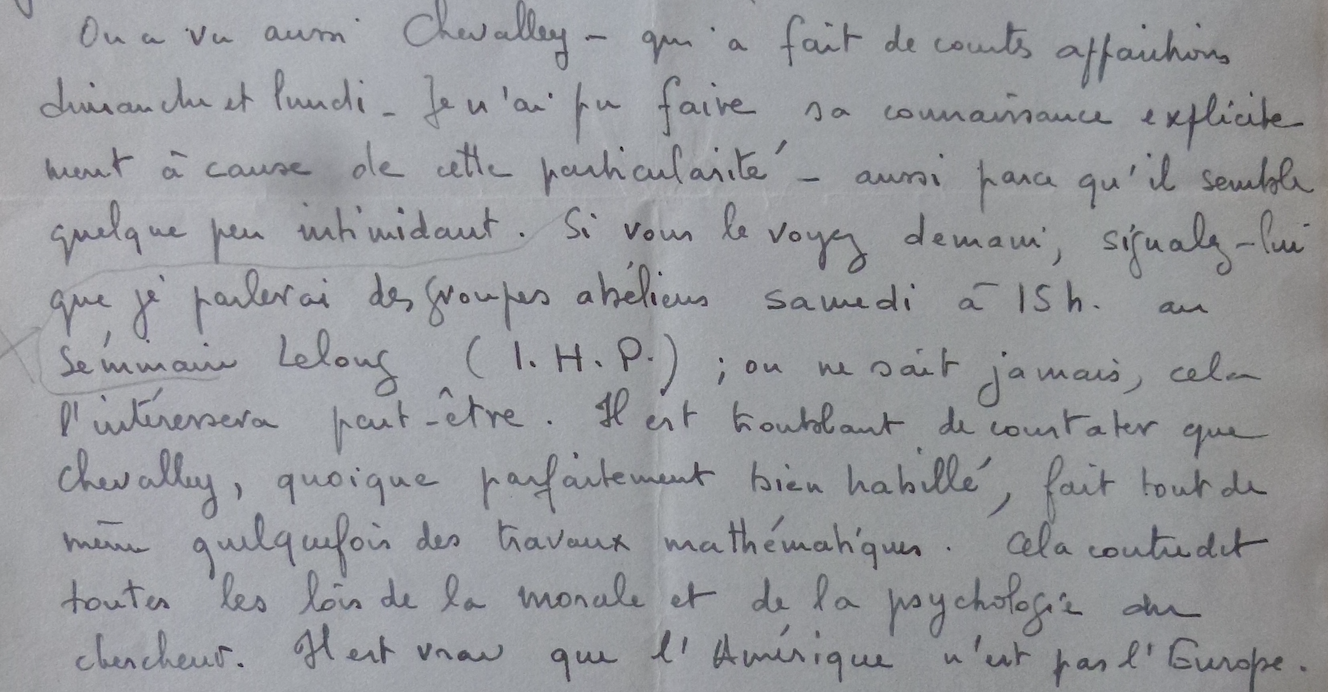}
\includegraphics[width=0.7\linewidth]{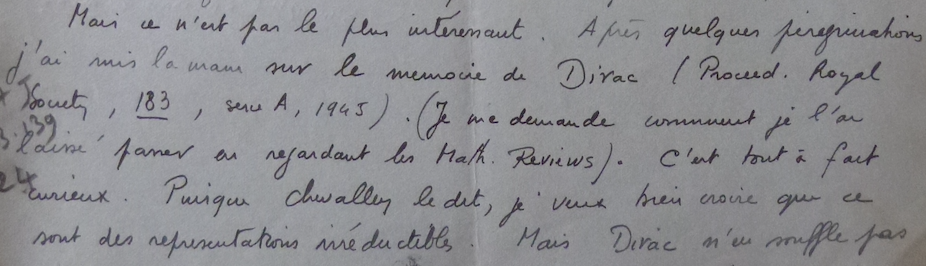}\end{center}
\caption{Rencontre avec Chevalley (mai 1946)}
\end{figure}

Le même Chevalley, alors en poste à Princeton, vient à Paris deux mois plus tard. C'est visiblement lui qui informe Godement des travaux tout récents de Dirac~\parencite*{Dirac_45}, alors en visite à Princeton, sur les représentations du groupe de Lorentz $\mathrm{SO}(3,1)$. Ces travaux, et ceux d'Harish-Chandra qui fait alors sa thèse avec Dirac sur un sujet proche, auront une influence profonde sur l'évolution ultérieure de la théorie des représentations. 

À Princeton toujours, le physicien Valentine Bargmann publie en~1947 une étude complète des représentations unitaires irréductibles du groupe $G = \mathrm{SL}(2, \R)$ ; il découvre notamment une série de représentations unitaires ayant des coefficients matriciels diagonaux de carré intégrable sur~$G$, et prouve des relations d'orthogonalité entre ces coefficients analogues aux relations de Schur pour les représentations des groupes finis \parencite{Bargmann}. Sitôt mis au courant, Godement établit des relations d'orthogonalité du même type pour les coefficients matriciels de représentations de carré intégrable sur un groupe localement compact quelconque ; il utilise pour cela un beau résultat de sa thèse  sur les fonctions de type positif de carré intégrable~\parencite[th.~17]{Godement_these}.

\begin{figure}[h]
\begin{center}\includegraphics[width=0.8\linewidth]{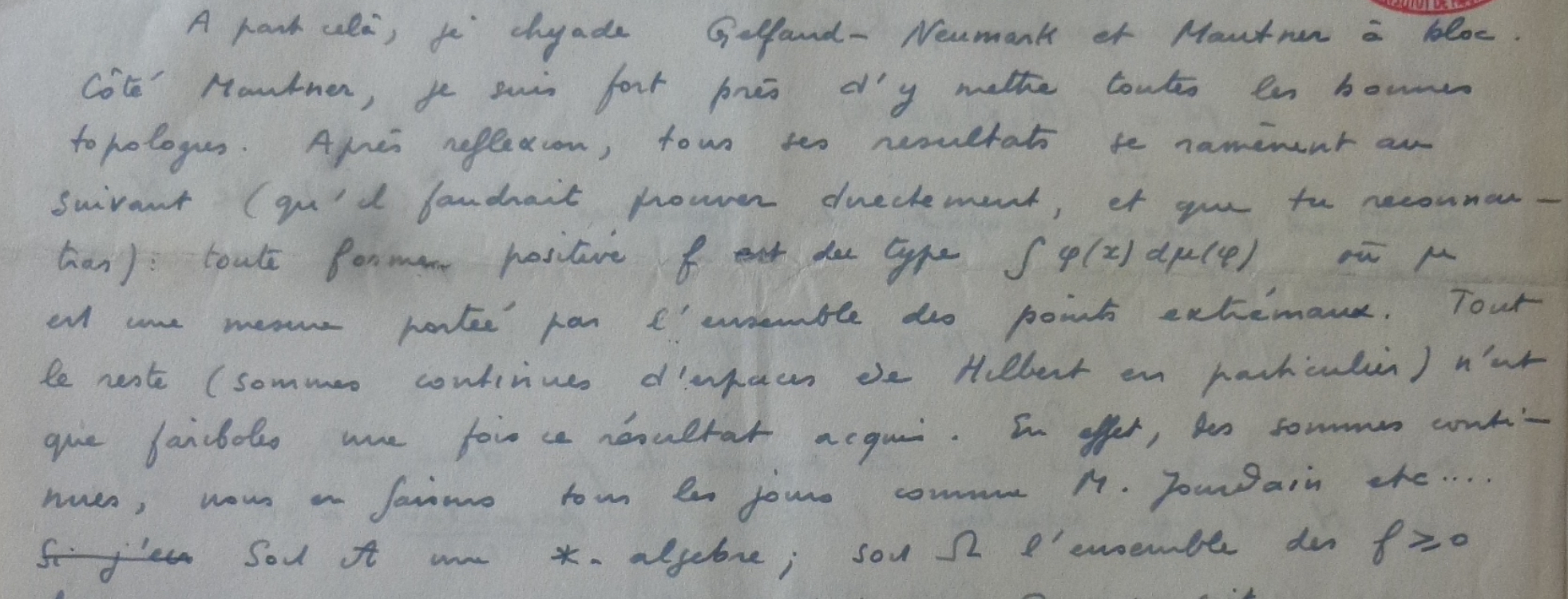}\end{center}
\caption{Tout près de Plancherel (novembre 1948)}\label{plancherel}
\end{figure}
La correspondance donne ainsi l'impression que, s'il est vrai que les idées mathématiques traversent l'Atlantique plus facilement que le rideau de fer, cela se fait à l'occasion des voyages des uns et des autres. De même, Cartan se rend aux États-Unis en~1948 ; à Harvard,  il rencontre George Mackey et écoute Friedrich Mautner. Les lettres de Godement montrent un intérêt mêlé de méfiance pour les travaux de Mautner~\parencite*{Mautner} ---  il lui semble que des méthodes de type Krein--Milman pourraient remplacer avantageusement ses arguments d'algèbres d'opérateurs. On sent la même distance aux travaux d'Irving Segal, qui songeait à utiliser des algèbres d'opérateurs en théorie des groupes dès 1940~\parencite{Segal_41} mais  ne publie la plupart de ses résultats qu'après~1947.

Un an plus tard, c'est Mackey  qui passe plusieurs mois à Nancy et informe Godement des résultats de Segal et Mautner sur la décomposition des représentations («\,formule de Plancherel\,» : \cite{Segal_Plancherel}), que Godement avait pressentis sinon démontrés. 

\begin{figure}[h]
\begin{center} \includegraphics[width=0.65\linewidth]{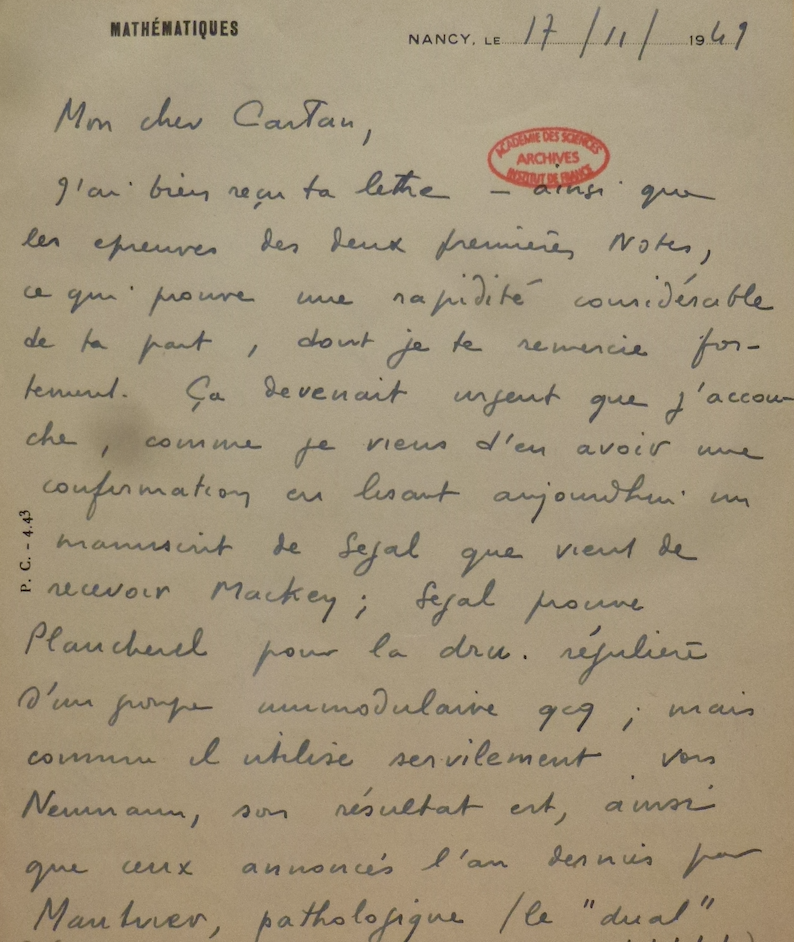}\end{center}
\caption{Concurrence avec Segal (novembre 1949)}
\end{figure}

Ce mélange d'échanges d'idées et de concurrence ne durera plus beaucoup : les derniers textes de Godement sur la théorie générale des représentations unitaires (consacrés à la notion de caractère pour les représentations de dimension infinie\footnote{Rappelons que si $\pi\colon G \to \mathcal{U}(H)$ est une représentation unitaire de~$G$, et si  l'espace~$H$ est \emph{de dimension finie}, alors le caractère de~$\pi$ est la fonction $g \mapsto \mathrm{Tr}(\pi(g))$, de~$G$ dans~$\C$ ; il résulte du th.~\ref{th:realisation_hilbertienne} qu'elle est de type positif. Si~$H$ est de dimension infinie, les opérateurs~$\pi(g)$ n'ont pas de trace en général ; la formule précédente ne permet donc pas de définir une fonction jouant le rôle du caractère de~$\pi$. De 1948 à 1954, Godement étudie certaines mesures « de type positif » sur~$G$, ou certaines distributions sur~$G$ dans le cas où~$G$ est un groupe de~Lie, pouvant jouer un rôle analogue à celui des caractères pour les représentations unitaires de dimension infinie. Si~$G$ est un groupe de~Lie semi-simple, Harish-Chandra adopte en 1951 une définition des caractères des représentations irréductibles qui utilise la théorie des distributions de Schwartz  (indépendamment de Godement, mais suite à des discussions avec Chevalley). Les travaux d'Harish-Chandra sur ces caractères, qui s'étalent sur près de trente ans, restent vu d'aujourd'hui l'un des sommets du sujet.  }) paraissent en 1954, juste avant son retour à Paris. Il cesse vite de travailler sur les représentations à ce niveau de généralité («\,abstract nonsense\,», dira-t-il) et s'intéresse aux formes automorphes : pour tout cela, voir l'entretien paru dans la \emph{Gazette}~\parencite{Entretien}. Peu après 1950, ses efforts de rédaction se concentrent sur un livre sur les faisceaux~\parencite*{Godement_faisceaux}, qui rendra grand service ; et sur Bourbaki (livre d'Intégration, exposés au séminaire). 

\begin{figure}[h]
\begin{center} \includegraphics[width=0.7\linewidth]{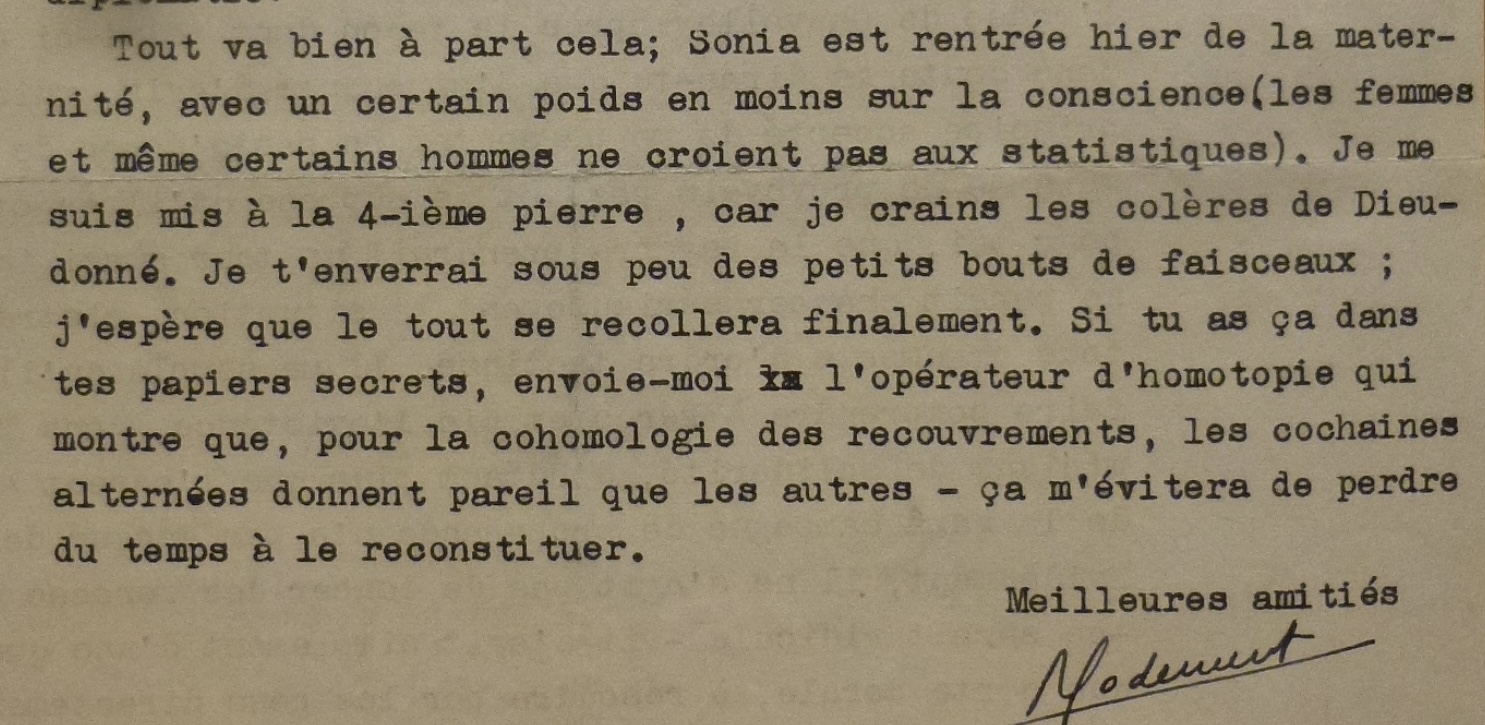}\end{center}
\caption{Famille et faisceaux (juin 1954)}
\end{figure}

\section{Et aujourd'hui ?}\label{sec:stochastiques}

Aujourd'hui, il me semble qu'on ne peut pas dire que les fonctions de type positif soient beaucoup étudiées pour elles-mêmes. Les généralités à leur propos se rencontrent surtout quand on s'instruit sur la théorie générale des représentations, sur les $C^\ast$-algèbres, ou sur les processus stochastiques ; en revanche, certaines fonctions de type positif jouent un rôle essentiel dans divers domaines. 

Pour ce qui est des représentations, les théorèmes~\ref{th:realisation_hilbertienne} à~\ref{th:gelfand_raikov} restent peut-être les résultats les plus utiles portant sur tous les groupes localement compacts. Ils servent surtout aux fondements de la théorie. Cela dit, beaucoup de travaux sur des classes particulières de groupes nécessitent une étude fine de \emph{certains} coefficients matriciels diagonaux de représentations ; ces fonctions sont de type positif comme on l'a vu, mais les problèmes qui se présentent alors dépendent beaucoup de la structure du groupe étudié. 

C'est le cas, par exemple, dans les travaux d'Harish-Chandra ou de Langlands concernant les représentations des groupes de Lie réductifs : l'importance de ces travaux (notamment pour l'étude des formes automorphes) n'est plus à démontrer, et le comportement de certains coefficients matriciels diagonaux y joue un rôle crucial --- mais les méthodes les plus efficaces pour aborder ce sujet ne se rattachent pas vraiment aux généralités sur les fonctions de type positif. 

D'autre part, beaucoup de «\,fonctions spéciales\,» de la physique mathématique (fonctions de Bessel, Jacobi, Legendre, et tant d'autres...) apparaissent comme coefficients matriciels diagonaux de représentations irréductibles. C'est fort utile pour comprendre et étudier les propriétés de ces fonctions à l'aide de la structure précise des groupes correspondants : un bon exemple est le texte encyclopédique de Vilenkin et Klimyk~\parencite*{Vilenkin}.

L'actualité récente concernant les fonctions de type positif sur des classes générales de groupes est probablement à chercher du côté des algèbres d'opérateurs. 

Par exemple, les progrès récents sur les algèbres de von~Neumann ont permis de mieux comprendre les fonctions de type positif sur les réseaux des groupes de Lie semi-simples. Or, la fonction indicatrice d'un sous-groupe discret~$\Gamma$ d'un tel groupe~$G$ est une fonction de type positif sur~$\Gamma$,  et cette fonction est constante sur les classes de conjugaison si et seulement si le sous-groupe~$\Gamma$ est distingué. Par conséquent, l'étude générale des fonctions de type positif invariantes par conjugaison sur~$G$ fournit des informations sur la structure des sous-groupes discrets distingués~$\Gamma$ de~$G$, sur les actions ergodiques de tels~$\Gamma$ sur des espaces probabilisés, sur les sous-groupes discrets distingués aléatoires de~$G$ (travaux de Bader, Bekka, Boutonnet, Gelander, Houdayer,  Peterson...). Voir l'exposé de Cyril Houdayer~\parencite*{Houdayer} au dernier Congrès international. 

Pour donner un second exemple, considérons la variante suivante de la notion de fonction de type positif sur un groupe localement compact~$G$ : disons qu'une fonction~$f\colon G \to \C$ est un \emph{multiplicateur complètement borné} sur~$G$ s'il existe un espace hilbertien~$H$ et des applications continues bornées $\varphi, \psi\colon G \to H$ telles que l'on ait $f(g^{-1}h)=\langle \varphi(h), \psi(g)\rangle$ pour tout $(g, h) \in G^2$.  Depuis un travail de Cowling et Haagerup~\parencite*{Cowling_Haagerup}, ces fonctions ont pris une importance croissante en théorie (géométrique) des groupes et en algèbres d'opérateurs, reliée à l'essor des divers visages de la propriété~$]\textup{T}[$ de Kazhdan pour les groupes localement compacts. Elles ouvrent la voie à la notion de \emph{moyennabilité faible} pour les groupes et les algèbres d'opérateurs  (de~Cannière, Cowling, Haagerup, Pisier, Ozawa, de~Laat, Knudby...) ; voir par exemple la synthèse récente d'Ignacio Vergara~\parencite*{Vergara}.

\smallskip

Terminons en évoquant les processus stochastiques, car si vous aviez déjà rencontré les fonctions de type positif sans que ce soit en lien avec les représentations de groupes ou les algèbres d'opérateurs, c'était peut-être comme \emph{fonctions de covariance de champs aléatoires}. 

Convenons d'appeler \emph{champ aléatoire} sur un ensemble~$E$ toute variable aléatoire dont les valeurs sont des fonctions de $E$ dans~$\R$ (on disait naguère « fonction aléatoire »). Si~$\Phi$ est un champ aléatoire sur~$E$, alors pour tout~$x \in E$, la valeur $\Phi(x)$  définit une variable aléatoire à valeurs dans~$\R$. Supposons de plus  que toutes les variances $\mathbb{V}[\Phi(x)]$, $x \in E$, soient \emph{finies}. Sans perte de généralité, on peut supposer qu'il existe un espace probabilisé~$\Omega$  tel que toutes les variables aléatoires $\Phi(x)$ s'identifient à des fonctions mesurables de $\Omega$ dans~$\R$ ; l'hypothèse est alors que $\Phi(x) \in L^2(\Omega; \R)$ pour tout $x \in E$. 

La \emph{fonction de covariance de~$\Phi$} est la fonction  $C\colon E\times E \to \R$ définie par $C(x,y) = \mathbb{E}[\Phi(x) \Phi(y)] = \langle \Phi(x), \Phi(y)\rangle_{L^2(\Omega;\R)}$. Elle ne dépend que de la \emph{loi} du processus~$\Phi$, et elle est \emph{toujours de type positif} : en effet, pour tous $n \geq 1$ et  $(x_1, \dots, x_n) \in E^n$, la matrice  $(C(x_i,x_j))_{1 \leq i,j \leq n}$ est égale à $\mathrm{Gram}_{L^2(\Omega;\R)}(\Phi(x_1), \dots, \Phi(x_n))$. Le fait que $C$ soit toujours de type positif était connu déjà de Michel Loève~\parencite*{Loeve}, au moins quand $E = \R$ ; et peut-être n'était-il pas le premier. 

C'est ainsi que les fonctions de type positif interviennent dans l'étude de processus stochastiques parmi les plus célèbres (y compris le mouvement brownien). Par exemple, si~$E$ est un espace topologique compact et si la fonction~$C$ est continue sur $E \times E$, alors les résultats de Mercer (1909) impliquent aussitôt l'existence d'un \emph{développement de Karhunen--Loève} pour le processus~$\Phi$ (voir par exemple~\cite[ch.~3]{AdlerTaylor}). De tels développements ont une portée pratique considérable de nos jours ; mais c'est une autre histoire.

\printbibliography % ou 
%\printauthorsdetails 

%\bibliographystyle{plainnat}
\flushleft

\end{document}